\documentclass[english,12pt,leqno]{article}
\usepackage[latin1]{inputenc}

\usepackage{fullpage}

\usepackage{tipa}
\usepackage{dsfont}
\usepackage{tikz-cd}
\usepackage{float}
\usepackage{graphicx}

\usepackage{xr}
 
\usepackage{amsxtra}
\usepackage{amsmath}
\usepackage{amssymb}
\usepackage{amsfonts}
\usepackage[all,2cell]{xy}
\usepackage{mathrsfs}
\usepackage{amsthm}
\usepackage{enumitem}
\usepackage{hyperref}
\usepackage{subfigure}
\usepackage{mathabx}
\usepackage{euscript}
\usepackage[bbgreekl]{mathbbol}
\usepackage{pgfplots}

\usepgfplotslibrary{fillbetween}

\usetikzlibrary{patterns}
\usetikzlibrary{calc}
 
\makeatletter
\hypersetup{
    unicode=true,           
    pdftoolbar=true,        
    pdfmenubar=true,        
    pdffitwindow=false,     
    pdfstartview={FitH},    
    pdftitle={\@title},     
    pdfauthor={\@author},   
    pdfsubject={},          
    pdfcreator={},          
    pdfproducer={},         
    pdfkeywords={},         
    pdfnewwindow=true,      
    colorlinks,             
    linkcolor=black,        
    citecolor=black,        
    filecolor=black,        
    urlcolor=black          
}
\makeatother

\newtheorem{thm}[equation]{Theorem}
\newtheorem{cor}[equation]{Corollary}
\newtheorem{lem}[equation]{Lemma}
\newtheorem{prop}[equation]{Proposition}

\newtheoremstyle{example}{\topsep}{\topsep}%
     {}
     {}
     {\bfseries}
     {.}
     {2pt}
     {\thmname{#1}\thmnumber{ #2}\thmnote{ #3}}

   \theoremstyle{example}
   
   \newtheorem{Defi}[equation]{Definition}
   \newtheorem{rem}[equation]{Remark}

   \newtheorem{exas}[equation]{Examples}
   \newtheorem{ex}[equation]{Example}

\numberwithin{equation}{section}

\def\AAA{\mathbb{A}}

\def\CC{\mathbb{C}}

\def\PP{\mathbb{P}}
\def\RR{\mathbb{R}}
\def\ZZ{\mathbb{Z}}
\def\QQ{\mathbb{Q}}

\def\gen{\mathfrak{g}}

\def\Aen{\mathfrak{A}}

\def\Bc{\mathcal{B}}

\def\Fc{\mathcal{F}}

\def\Lc{\mathcal{L}}

\def\Hc{\mathcal{H}}
\def\Oc{\mathcal{O}}

  \def\Coker{\operatorname{Coker}\nolimits}
 
 \def\conf{{\mathfrak{conf}}}
  \def\Conf{{\on{Conf}}}

 \def\dbar{{\ol{\partial}}}
 \def\del {{\partial}}

 \def\ent{{\on{ent}}}

 \def\hol{{\on{hol}}}
  \def\Hom{\operatorname{Hom}\nolimits}

\def\is{{\on{is}}}

 \def\kap{{\varkappa}}
  \def\Ker{\operatorname{Ker}\nolimits}

 \def\Lie{{\on{Lie}}}
   \def\lra{\longrightarrow}

\def\on{\operatorname}
\def\ol{\overline}

 \def\so{{\mathfrak{so}}}
 \def\Spec{{\on{Spec}}}

 \def\Tot{{\on{Tot}}}
 
\def\wt{\widetilde}

\def\1{{\mathbf{1}}}
\def\2{{\mathbf{2}}}
\def\(({(\hskip -1mm (}
\def\)){)\hskip -1mm )}
 \def\be{\begin{equation}}
\def\ee{\end{equation}}
\def\ed{\end{document}}

\def\-{{\setminus}}

\title{ Conformal maps in higher dimensions and derived geometry}
\author{  Mikhail Kapranov }
 
\begin{document}

\maketitle

\begin{abstract}
 By Liouville's theorem,  in dimensions $3$ or more conformal
 transformations form a finite-dimensional group, an apparent drastic departure
 from the $2$-dimensional case. We propose a derived enhancement of
 the conformal Lie algebra which  is an infinite-dimensional
  dg-Lie algebra incorporating not only symmetries but also
 deformations of the conformal structure.
   Our approach is based on (derived) deformation
 theory of the ambitwistor  space of complex null-geodesics. 
\end{abstract}
 

\section*{Introduction}

The classical theorem of Liouville (1850) says that the behavior of conformal maps
in dimensions $\geq 3$, as compared to dimension $2$, is drastically different.
While in $2$ dimensions, holomorphic functions give an infinite-dimensional supply of
local conformal maps, any conformal map between connected domains $U,V\subset\RR^n$,
$n\geq 3$, comes from a global M\"obius transformation of the conformal sphere $S^n$,
i.e., from an element of the finite-dimensional group $O(n+1,1)$. This
apparent discontinuity can be rather puzzling.

\vskip .2cm

The goal of this note is to propose a resolution to this apparent puzzle by using
the point of view of derived geometry, i.e., of the homological, derived category-style
approach to algebraic and differential geometry 
\cite{drinfeld-letter, HS-def, lurie-DAG, TV}.  More precisely, we recover the 
missing infninite-dimensional part of the conformal group (it is technically easier 
to  start with the Lie algebra) in a  different cohomological degree. For $n=2$ this
degree is $0$, and we observe infinite-dimensionality at the naive classical level. 

\vskip .2cm

For an $n$-dimensional
 complex analytic conformal manifold $M$ of appropriate type we introduce a differential
 graded (dg-) Lie agebra $R\conf(M)$ which is of infinite-dimensional nature regardless
 of $n$. Its cohomology includes:
 
 \begin{itemize}
 \item $H^0=\conf(M)$, the usual Lie algebra of conformal Killing vector fields 
 (locally, infinite-dimensional
 for $n=2$, finite-dimensional for $n>2$).
 
 \item $H^1$ being the space of infinitesimal deformations of the conformal structure
 (locally, zero for $n=2$, infinite-dimensiomal for $n>2$). 
 \end{itemize}
 
 Thus, for $M=\CC^n$ with flat metric, the total size of $H^\bullet R\conf(M)$ varies
 ``continuously'' with $n$. To illustrate this point, we identify $H^\bullet R\conf(\CC^n)$
 as a representation of $SO(n,\CC)$ in Theorem \ref{thm:conf-flat}. 

\vskip .2cm

Our approach is based on the analogy with  behavior of holomorphic functions on
$\CC^n$ vs. $\CC^n\setminus \{0\}$. While  passing from $\CC$ to $\CC\setminus\{0\}$
increases the supply of holomorphic functions, passing from $\CC^n$ 
to $\CC^n\setminus\{0\}$,  $n\geq 2$,  does not (Hartogs' theorem). But if we look at the
total cohomology $H^\bullet(\CC^n\setminus\{0\}, \Oc)$ of the sheaf of
holomorpic functions, we find the missing singular parts in the cohomological degree
$n-1$. 

\vskip .2cm

To relate Liouvile's theorem with Hartogs-type phenomena, we   use the ambitwistor
approach to holomorphic conformal geometry \cite{LB-tams, manin, mason}. The principal object there
is $L(M)$, the space of complex null-geodesics in a holomorphic conformal manifold $M$,
with its natural contact structure. Its holomorphic contact geometry completely encodes
the holomorphic conformal geometry of $M$. 
In particular, the space of (holomorphic) conformal
Killing fields on $M$ is found as the space of sections
\be\label{eq:H0=conf}
\conf(M) \, = \, H^0 (L(M), \kap), 
\ee
where $\kap$ is the sheaf of holomorphic contact vector fields on $L(M)$. 
Already in the local case ($M$ is a small geodesically convex neighborhood of a point $x_0$),
$L(M)$ is a complex manifold which (for $n\geq 3$) is not compact but has compact
directions: it contains $(n-2)$-dimensional complex projective quadrics $L_x$
formed by null-geodesics passing through various points $x\in M$. For
manifolds of this type,
 coherent sheaves such as $\kap$, can have  finite-dimensional $H^0$ but
infinite-dimensional higher
cohomology. We define
\[
R\conf(M) \,=\, R\Gamma(L(M), \kap),
\]
the dg-Lie algebra of derived global sections of the sheaf of Lie algebras $\kap$.

\vskip .2cm

There can be other approaches to defining $R\conf(M)$, for example, in the $C^\infty$,
rather than holomorphic case. The holomorphic approach adopted here provides a natural
way to arrive at the idea of a derived extension of the conformal algebra. It also
leads to a conceptually transparent proof of  (the infinitesimal, holomorphic version of)
Liouvuille's theorem, in the form of finite-dimensionality of the $H^0$-space
in \eqref{eq:H0=conf}.

\vskip .2cm

Since for $n\geq 3$, the infinite-dimensionality of the dg-Lie algebra $R\conf(M)$ is situated in the
odd cohomological degree $1$, one can integrate it, in a purely algebraic way, cf. 
\cite{calaque, lurie-moduli}, 
to a {\em derived group} $R\Conf(M)$. The classical truncation' of $R\Conf(M)$
is $\Conf(M)$, the usual Lie group of (holomorphic) conformal diffeomorphisms, and the
whole $R\Conf(M)$ can be seen as an infinite-dimensional formal derived thickening
of $\Conf(M)$.  Such derived groups, and their analogs for conformal superspaces 
\cite{deligne-freed, manin}
 should  be of importance for the study of (super)conformal quantum field theories in
dimensions $n\geq 3$. See \cite{saberi}
 for a  somewhat different recent appearance of derived geometry
constructions in that context. 

\vskip .2cm

I am grateful to C. Schweigert and M. Yamazaki for interest in this work and useful suggestions. 
This research was supported by the  World Premier International Research Center Initiative (WPI Initiative),   MEXT, Japan. 


\section{ Theorems of Hartogs and Liouville}

\paragraph{The classical theorems.}

The classical theorem of Hartogs says:

\begin{thm}
For $n \geq 2$, every holomorphic function on $\CC^n - \{0\}$ extends holomorphically to $\CC^n$. 
\end{thm}

In the algebro-geometric version, over the base field $\CC$ and over the Zariski topology,
we consider the affine space $\AAA^n$ (i.e., ``$\CC^n$ `considered as an algebraic variety''). 
The corresponding statement is that
for $n\geq 2$
\be\label{eq:hartogs}
H^0 (\AAA^n\-  \{0\}, \Oc) \, = \,  H^0 (\AAA^n, \Oc)\, =\,  \CC[z_1,\cdots, z_n],
\ee
  (no increase of the ring of regular functions), while for $n=1$ we have a manifest increase:
  \be\label{eq:mer:m=1}
  \CC[z, z^{-1}] \,=\,  H^0(\AAA^1 \-  \{0\}, \Oc) \, \supsetneq \,  H^0 (\AAA^1, \Oc) \,=\, \CC[z].
  \ee
  with tne new  part being $z^{-1}\CC[z^{-1}]$ (space of polar parts of functions with poles at $0$). 
  
  This phenomenon looks like discontinuity:
   something seemingly disappears as we pass to higher dimensions. 
  
  \vskip .2cm
  
  There is an even more classical result in geometry where the situation changes drastically
  in passing to higher dimensions: the Liouville theorem. We consider the flat
  Euclidean space   $(\RR^n, \langle -,-\rangle)$ and look at  conformal transformations between
  domains $U$, $V$ in $\RR^n$. 
  If  $n=2$, then locally, we have an infinite-dimensional supply of such transformations, as
  any holomorphic function on a domain in $\CC=\RR^2$   is conformal. However, the Liouville theorem says:
  
  \begin{thm}
  Let $n\geq 3$, let  $U,V$ be connected open domains in $\RR^n$  and  $\varphi: U\to V$ be a conformal
  diffeomorphism. Then $\varphi$ extends to a  Moebius-type conformal map (composition of rigid motions, dilations and inversions)
   defined on the entire $\RR^n$
  with the possible exception of one point, and belonging to the standard conformal group
  $O(n+1,1)$. 
  \end{thm}

  Let us recall the geometric meaning of the group $O(n+1,1)$ in this case.  For this, and for further
  analysis, it is convenient to work in the complex analytic situation.
  
 \paragraph{The complex setting.} 
  
  Let $M$ be a complex manifold of   dimension $n$.  We can speak about {\em holomorphic Riemannian metrics} on $M$.
  Such a metric is a holomorphic section of the vector bundle $S^2(T^*_M)$ which is
  non-degenerate at each point. In holomorphic coordinates it is gives as a symmetric matrix
  $g(z)=\|g_{ij}(z)\|$ of holomorphic functions. 
  
  A {\em (holomorphic) conformal metric} on $M$ is, naively,  a ``Riemannian metric
  defined up to a scalar''.  This means that locally we have a representative $g_{ij}(z)$ which is
  a holomorphic Riemanian metric, with the understanding that we identify
  \[
  g_{ij}(z) \,\sim \, \lambda(z) g_{ij}(z), \quad \lambda \in \Oc_M^*,
  \]
  i.e., that $ \lambda(z) g_{ij}(z)$ represents the same conformal metric  
  as $g_{ij}(z)$.

 If we fix a point $z\in M$, then a complete  invariant of a ``non-degenerate symmetric form on $T_zM$ defined modulo scalars'',
 is the {\em null-cone}  $C_zM \subset T_zM$.  Thus, globally:
 
 \begin{Defi} (a)  A {\em holomorphic conformal metric}
 on a complex manifold $M$ is a holomorphic family of non-degenerate quadratic cones $C_z\subset T_zM$. 
 A {\em complex conformal manifold} is a complex manifold equipped with a holomorphic conformal metric.
 
 \vskip .2cm
 
 (b) Let  $(M, (C_z))$ and $(M', (C'_{z'}))$  are two  $n$-dimensional complex  conformal manifolds  
 A {\em conformal mapping} $\varphi: M\to M'$ is a biholomorphic map whose
 differential takes null-cones to null-cones: $d_z\varphi(C_z) = C_{\varphi(z)}$. 
 \end{Defi}

 \vskip .2cm
  
  An alternative definition would be that  a conformal metric is given by a holomorphic line bundle 
  $\Lambda$ on $M$  
  an $\Lambda$-valued scalar product  $g\in \Hom(S^2T_M, \Lambda)$ on $T_M$. 
  It is easily seen to be equivalent to the above. 
  
  \begin{ex}[(Conformal quadric)] 
(a)  Let $Q\subset \PP^{n+1} = \PP(\CC^{n+2})$ be a non-degenerate quadric hypersurface. 
 It has a canonical conformal structure defined as follows.  For $z\in \QQ$ let $\ol T_z Q\subset \PP^{n+1}$ be the
 {\em projective tangent space} to $Q$ at $z$. It is an algebraic variety  isomorphic to $\PP^n$,
 and we have a canonical identification of the (usual) tangent spaces $T_z \ol T_z Q = T_z Q$. 
 The intersection $\ol C_z = Q\cap \ol T_zQ$  (inside $\PP^{n+1}$) 
 is a quadratic hypersurface in $\ol T_zQ$ with
 one singular point, namely $z$. It is called the {\em projective tangent cone} to $Q$. 
 The (intrinsic) tangent cone to this intersection is a nondegenerate
 quadratic cone $C_z$ in $T_z \ol T_z Q$, i.e., in $T_zQ$. 
 This defines the conformal structure on $Q$. 
 
 \vskip .2cm
 
 (b) Thus any automorphism of $\PP^{n+1}$ preserving $Q$ defines a conformal mapping $Q\to\QQ$.
 Such automorphisms form the group $O(n+2,\CC)$. 
 
 \vskip .2cm
 
 (c) If we fix one point $\infty\in Q$, then $Q \- \ol C_\infty Q$ is isomorphic to $\AAA^n$
 as an algebraic variety so to $\CC^n$ as a complex manifold. 
 The induced conformal structure on $\CC^n$ is the standard flat conformal structure.
 Thus any $g\in O(n+2,\CC)$ defines a conformal mapping from $\CC^n$ (minus, possibly,
 a quadratic cone hypersurface) to itself.
 
 \vskip .2cm
 
 (d) If we want to restrict to real points, then $Q$ gives the $n$-sphere $S^n$, with its
 standard conformal structure. The real group $O(n+1,1)$ acts by conformal
 mappings of $S^n$ to itself. 
 Each  projective tangent cone $\ol C_zS^n$
 gives just the point $z$, and $S^n -  \{\infty\}$ is identified with $\RR^n$ via
 the stereographic projection. In this way any $g\in O(n+1,1)$ defines
 a conformal map of $\RR^n$ (minus, possibly, a single point) to itself. 
 \end{ex}
 
 The holomorphic version of Liouville's theorem can be formulated as follows.
 
 \begin{thm}
 Let $n\geq 3$, let $U,V\subset Q$ be open domains and $\varphi: U\to V$ be
 a holomorphic conformal mapping. Then $\varphi$ extends to an isomorphism
 of algebraic varieties $\Phi: Q\to Q$. 
 \end{thm}
 
 Let us concentrate on the Lie algebra (infinitesimal) version. Let $M$ be a complex manifold with
 conformal structure. A {\em conformal Killing field} on $M$ is a holomorphic vector field
 preserving the conformal structure. 
 We denote by $\conf(M)$ the Lie algebra formed by conformal Killing fields. 
 The Lie algebra version of the Liouville theorem is:
 
 \begin{thm}\label{thm:liou-inf}
 Let $Q\subset \PP^{n+1}$ be the $n$-dimensional quadric.
 
 \vskip .2cm
 
 (a) For any $n\geq 1$ we have $\conf(Q)= \so(n+2,\CC)$.
 
 \vskip .2cm
  (b) Let $n\geq 3$. 
 For any connected open domain $U\subset Q$ the restriction map 
 \[
 \so(n+2,\CC)=\conf(Q) \lra \conf(U)
 \]
  is
 an isomorphism. 
 \end{thm}
 
 As we shall see, this fact is a manifestation of the same phenomenon as the Hartogs theorem and can
 be overcome in a similar way. 
 
  
  \section{Overcoming Hartogs}
  
  A way to recover the missing polar parts  in Hartogs' theorem is by using the full cohomology $H^\bullet (\AAA^n\setminus \{0\}, \Oc)$, not just $H^0$.
 We have the following elementary fact.
 
 \begin{prop}\label{prop:An-0-O}
 For any $n\geq 2$ we have
 \[
 H^i(\AAA^n \setminus \{0\}, \Oc) \,\simeq\,  \begin{cases}
 \CC[z_1,\cdots, z_n], & \text{ if } i=0, 
 \\
 z_1^{-1}\cdots z_n^{-1} \cdot \CC[z_1^{-1}, \cdots, z_n^{-1}], & \text { if } i=n-1, 
 \\
 0, & \text{ otherwise}. 
 \end{cases}
 \]
 \end{prop}
 
 Here the new space, formed by  the polar parts
 \be\label{eq:H-n-0}
 H^{n-1}(\AAA^n\-\{0\}, \Oc) \,=\, H^n_{\{0\}} (\AAA^n, \Oc)\,=\, z_1^{-1}\cdots z_n^{-1} \cdot \CC[z_1^{-1}, \cdots, z_n^{-1}]
 \ee
 appears in cohomological degree $n-1$ and is therefore invisible if we remain in the classical (non-homological)
 framework. It can also be seen as the $n$th cohomology with support at $0$. In fact, the second identification
 in \eqref{eq:H-n-0} holds for any $n \geq 1$. Thus passing to the cohomology restores the continuity.
 
 \vskip .2cm
 
 The easiest way to establish Proposition \ref{prop:An-0-O} is by using the  \v{C}ech complex associated to the covering
 of  $\AAA^n\-\{0\}$ by the affine open sets $U_i = \{z_i\neq 0\}$, $i=1,\cdots, n$. The space
 $\Gamma( U_{i_1, \cdots, i_p},\Oc)$
 of regular functions on
 each $p$-fold  intersection  is realized inside the Laurent polynomial ring 
 $\CC[z_1^{\pm 1},\cdots, z_n^{\pm 1}]$, and the $(n-1)$th cohomology will appear as the span of the Laurent
 monomials which will not appear in any  of the   $\Gamma( U_{i_1, \cdots, i_p},\Oc)$. This approach is equivalent
 to the classical computation of the cohomology of the sheaves $\Oc(d)$ on $\PP^{n-1}$ due to Serre,
 see \cite{hartshorne-AG}. 
 
 \vskip .2cm
 
 If we consider the complex manifiold $\CC^n$ instead of the algebraic variety $\AAA^n$ and the
 sheaf $\Oc_\hol$ of holomorphic functions, we have a statement similar to 
 Proposition \ref{prop:An-0-O}, but with
 \[
 H^{n-1}(\CC^n\-\{0\}, \Oc_\hol) \,=\, H^n_{\{0\}} (\CC^n, \Oc_\hol)
 \,=\, z_1^{-1}\cdots z_n^{-1} \cdot \CC[[z_1^{-1}, \cdots, z_n^{-1}]]_\ent
 \]
 being the space of Taylor series representing  entire functions. This  space is known as
 the space of holomorphic hyperfunctions on $\CC^n$ with support at $0$, see \cite{SKK}. 
 
 \vskip .2cm
 
 For any sheaf $\Fc$ of $\CC$-vector spaces on a topological space $X$ we denote by $R\Gamma(X,\Fc)$
 the derived functor of sections of $\Fc$, i.e., ``the'' complex of $\CC$-vector spaces whose cohomology is $H^\bullet(X.\Fc)$. 
 Such a complex is defined uniquely up to unique isomorphism in the derived category. 
 If $\Fc$ has some additional algebraic structure (commutative algebra, Lie algebra etc.),
  then  it is well known that $R\Gamma(X,\Fc)$ can be defined in such a way as
 to inherit this structure. 
 
 \begin{exas}\label{ex:algebra-str}
 (a)  Let $X$ is a complex manifold and $\Fc=\Oc_X$ be the sheaf of holomorphic functions.
 It is a sheaf of commutative algebras. 
 The Dolbeault complex $\Omega^{0,\bullet}(X), \dbar)$ is a model for $R\Gamma(X, \Oc_X)$ which
 has the structure of a commutative dg-algebra.
 
 The sheaf $\Fc=T_X$  of holomorphic vector fields on $X$ is a sheaf of Lie algebras.
 The Dolbeault complex $(\Omega^{0,\bullet}(X, T_X),\dbar)$ is a dg-Lie algebra model for
 $R\Gamma(X,T_X)$, with the Lie structure given by the Schouten bracket. 

\vskip .2cm

(b) Let $X$ be an algebraic variety with Zariski topology  and $\Fc=\Oc_X$ ibe the sheaf of regular functions, A commutative
dg-aglebra model for $R\Gamma(X, \Oc_X)$ can be obtained as  the goobal relative de Rham complex 
$\Gamma(J, \Omega^\bullet_{J/X})$
where $J\to X$ is a {\em Jouianolou torsor}, i.e., an affine variety which is made into a Zariski ocally trivial bundle
over $X$ with fibers being affine spaces and transition functions being affine transformations. See
\cite{BD} for a general discussion and \cite{FHK} for a concrete example with $X=\AAA^n\-\{0\}$. 

\vskip .2cm

(c)  If $X$ is  any topological space
 and $\Fc$ is any sheaf of commutative (resp. Lie, etc.) $\CC$-algebras,
then the \v{C}ech model for $R\Gamma(X,\Fc)$ produces  a cosimplicial commutative
(resp. Lie, etc.) $\CC$-algebra. There is a general  procedure of converting a cosimplicial algebra
of any given type into a dg-algebra of the same type using the Thom-Sullivan construction
involving polynomial differential forms on simplices. It provides the most general way
to make $R\Gamma(X,\Fc)$ to inherit the algebra structure present on $\Fc$. We refer to 
 \cite{HS-def} \S 5.2 for details. 
 \end{exas}
 
  Thus the correct  $n$-dimensional replacement of the algebra $\CC[z, z^{-1}]$ of Laurent
 polynomials is the commutative dg-algebra
 \[
 \Aen_{[n]} \,=\, R\Gamma(\AAA^n\-\{0\}, \Oc) 
 \]
 defined as in Example \ref{ex:algebra-str}(b). In particular, tensoring $\Aen_{[n]}$ with a finite-dimensional
 reductive Lie algebra $\gen$ leads to interesting higher-dimensional derived generalizations
 of Kac-Moody algebras \cite{FHK, {gwilliam-williams}}. 
 

\section{Ambitwistor description of conformal metrics}
\label{sec:ambi}

 We want to show that Liouville's theorem,  at least in its  complex, infinitesimal form \eqref {thm:liou-inf}, can be seen as a Hartogs-type phenomenon
 and therefore can be ``overcome'' by introducing cohomological degrees of freedom. 
 For this,  we recall the main points of the  ambitwistor approach 
 \cite{LB-tams, manin, mason}  to holomorphic conformal metrics
 (in any dimension, in particular  without assuming self-duality in dimension 4).

 \paragraph{The space of null-geodesics.} 
 
 Let $(M, g)$ be an $n$-dimensional complex manifold with a holomorphic Riemannian metric. We can then speak about
 {\em null-geodesics} in $M$ which are parametrized holomorphic curves $\gamma: U\to M$,
 $U\subset\CC$ open,  
 satisfying the complex version of the geodesic equation and such that $\gamma'(t)$ is isotropic everywhere.
 The  elementary but fundamental fact  is, see \cite{LB-tams} \S II.2:
 
 \begin{prop}\label{prop:geod-conf}
 For two conformally equivalent metrics $g(z)$ and $\lambda(z) g(z)$, 
  the null-geodesics are the same
 up to a re-parametrizaton. \qed
 \end{prop}
 
 Put differently, let $QTM\subset\PP(TM)$ be the quadric bundle formed by the null-directions
 in the tangent spaces $T_xM$, $x\in M$. It is a complex manifold of dimension $2n-2$. 
 The ``complex geodesic flow'' for $g(z)$ is the 1-dimensional complex foliation $\Lc$  on $QTM$
 whose leaves are the tangent lifts of null-geodesics for $g(z)$. Note that $QTM$ depends only
 on the conformal class of $g(z)$. Proposition \ref{prop:geod-conf} says that so does $\Lc$. 
 
 \vskip .2cm
 
 Let now $(M, (C_x))$ be a holomorphic conformal manifold of dimension $n$. We then have
 the quadric bundle $QTM\subset \PP(TM)$ with fibers $Q(T_xM)= \PP(C_x)\subset \PP(T_xM)$. 
 By the above, $QTM$ carries a canonical 1-dimensional holomorphic foliation $\Lc$
 whose leaves, are, locally, the lifts of complex null-geodesics for any holomorphic metric
 representing $(C_x)$. 
 
 \vskip .2cm
 
 The {\em space of null-geodesics} $L=L(M)$ is defined as the space of leaves of the foliation $\Lc$. 
 In the sequel we assume that this space of leaves exists, i.e.,  intuitively, the global
 behavior of complex null-geodesics in not too wild. More precisely, following \cite{LB-tams},
 we make the following
 
 \begin{Defi}\label{def:civil}
 A holomorphic conformal manifold $(M, (C_x))$ of dimension $n$ is called {\em civilized},
 if:
 \begin{itemize}
 \item[(1)] There is a Hausdorff complex manifold $L$ of dimension $2n-3$ and a holomorphic
 submersion $\rho: QTM\to L$ whose fibers are precisely the leaves of $\Lc$, with the
 property:
 
 \item[(2)] The restriction of $\rho$ to any quadric $Q(T_xM)$, $x\in M$, is a holomorphic embedding
 (that is, no complex null-geodesic passes through the same point twice). 
 \end{itemize}
 \end{Defi} 
 
 For a civilized $M$ the manifold $L=L(M)$ is defined uniquely up to a unique isomorphism.
 In the sequel we will assume that $M$ is civilized. Here are some examples.

 \begin{exas}\label{ex:civilized}
 
 (a) (Flat case, noncompact) $M=\CC^n$ with a flat conformal metric. In this case $L(\CC^n)$ is a closed subvariety in
 the space of all straight lines in $\CC^n$. It is  an algebraic variety which,
 for $n\geq 3$,  is not affine and not projective.  More precisely, let $Q^{n-2}\subset \PP^{n-1}=\PP(M)$ 
 be the quadric formed by null-lines in $M$ through $0$. Any null-line in $M$ can be seen, in
 a unique way, as a translation of a line passing through $0$. This means that $L(M)$
 is the total space of the vector bundle $l\mapsto M/l$ on $Q$, i.e., of the bundle whose
 fiber over the point $[l]$ represented by a null-line $l$ through $0$, is $M/l$. 
 The bundle $l\mapsto  M/l$ is in fact defined on the entire $\PP(M)=\PP^{n-1}$ and as such,
 is identified with
 $T_{\PP(M)}(-1))$,
  because of the ``Euler sequence''  \cite{GH}:
 \[
 0\to \Oc_{\PP(M)}(-1) \lra M\otimes \Oc_{\PP(M)} \lra T_{\PP(M)}(-1) \lra 0. 
 \]
(We recall that $\Oc_{\PP(M)}(-1)$ is the tautological line bundle, i.e., $l\mapsto l$ in the
 above notation.) So we conclude that
 \be\label{eq:L(Cn)}
 L(\CC^n) \,=\, \Tot\bigl(T_{\PP^{n-1}}(-1)|_{Q^{n-2}}\bigr). 
 \ee

 \vskip .2cm
 
 (b) (Flat case, compact) $M=Q^n\subset \PP^{n+1}$ is the $n$-dimensional projective quadric. In this
 case $L(Q^n)$  consists of all straight lines in $\PP^{n+1}$ which lie on $Q^n$. It is a projective algebraic variety
 identified with $G^{\on{is}}(2, \CC^{n+2})$, the Grassmannian of $2$-dimensional subspaces in $\CC^{n+2}$
 which are isotropic with respect to the quadratic form defining $Q^n$. 
 
 \vskip .2cm
 
 (c) (Local case) We can always replace $M$ by a sufficiently small neighborhood
  around a fixed point $x_0$
 (small with respect to the curvature data of the metric near $x_0$).  Then the situation will be similar to the
 flat case (a), so we get a civilized manifold. See \cite{LB-tams}, \S II.1. 
 \end{exas}
 
 For each $x\in M$ we define 
 \be\label{eq:Lx}
 L_x\,=\, \bigl\{ \gamma\in L \bigl| \,\, x\in\gamma\bigr\} \,
 \subset  \, L
 \ee
 to consist of null-geodesics that pass through $x$.
 The condition (2) of Definition \ref{def:civil} implies that $L_x$ is identified with
 the quadric $Q(T_xM)$, i.e., is isomorphic to  $Q^{n-2}$. The conformal geometry of $M$
 is encoded by the system of subvarieties $L_x$. That is, $x$ and $y$ are null-separated,
 if and only if $L_x\cap L_y\neq\emptyset$. Further, comparison with the flat case
 identifies the normal bundle of each $L_x$ in $L$. That is, with respect to any identification
 $L_x\simeq Q^{n-2}$ we have
 \be\label{eq:N-Lx-L}
 N_{L_x/L} \,\simeq \, T_{P^{n-1}}(-1)|_{Q^{n-2}}. 
 \ee
 The bundle in the RHS of \eqref{eq:N-Lx-L} is homogeneous (equivariant
 under automorphisms of $Q^{n-2}$), 
  so  one can write  \eqref{eq:N-Lx-L}  without reference to
 a particular way of identifying $L_x$ with $Q^{n-2}$. 
 
 
 \paragraph{$L(M)$ as a contact manifold.} Let $X$ be a complex manifold of odd dimension
 $2m+1$. We recall \cite{arnold, LB-tams, ovsienko} that a (holomorphic) {\em contact structure}
 on $X$ is a (holomorphic) vector subbundle $\Theta\subset T_X$ of rank $2m$ which is maximally
 non-integrable in the following sense. Let $\kap = T_X/\Theta$ be the quotient line bundle. 
 Then $\Theta$ is given by the vanishing of the tautological $\kap$-valued
 contact form $\theta: T_X\to\kap$.  A local trivialization of $\kap$ makes $\theta$ into a
 usual holomorphic $1$-form. The maximal non-integrability conditions means that
 \be
 \theta\wedge \underbrace{ d\theta \wedge \cdots \wedge d\theta}_{\text{$m$}} \,\neq \, 0 \,\text{ everywhere}.
 \ee
 This condition is known to be independent on the way we represent $\theta$ as a usual form
 by choosing a trivialization of $\kap$. More precisely, since $\theta$ is, intrinsically, a 1-form
 with values in a line bundle, $d\theta$ is not invariantly defined, but it is invariantly defined
 modulo  (the wedge ideal generated by) $\theta$. Therefore $\theta\wedge (d\theta)^{\wedge n}$ is invariantly defined
 as a  volume form with values in $\kap^{\otimes (m+1)}$,
 becuase it involves only $d\theta$ modulo $\theta$, as $\theta\wedge\theta=0$. 
 Since it is nowhere vanishing,  we
 get a canonical identification of line bundles
 \be
 \kap^{\otimes (m+1)} \,\simeq\,  \omega_X^{\otimes(-1)},
 \ee
 where $\omega_X$ is the line bundle of volume forms. 
 
 \vskip .2cm
 
  Another consequence of the same  remark  
 is that $d\theta$ is invariantly defined and non-degenerate on $\Ker(\theta)=\Theta$, that is,
 we have a non-degenerate skew-symmetric form
 \be
 d\theta: \Lambda^2(\Theta) \lra \kap. 
 \ee
If now $Z\subset X$ is a smooth hypersurface which is transversal to $\Theta$ everywhere,
then it carries the $1$-dimensional  {\em bicharacteristic foliation} $\Bc$,
with tangent spaces to the leaves being the $1$-dimensional subspaces
\[
\Bc_z \,=\, \Ker (d\theta|_{T_zZ \cap \Theta_z})  \,\subset \, T_zZ, \quad z\in Z.
\]
 It is classical (the ``contact reduction'') that the space of leaves of $\Bc$ (the space of bicharacteristics), if it exists, 
 is again a contact manifold, now of dimension $2m-1$. 
 
 \vskip .2cm
 
 Let now $(M, (C_x))$ is a holomorphic conformal manifold of dimension $n$. Then
 $T^*M$ is a symplectic manifold, so $\PP(T^*M)$ is a contact manifold \cite{arnold}
 which is identified with $\PP(TM)$ by the conformal structure. The hypersurface
 $QTM\subset \PP(TM)$ of null-directions is transversal to $\Theta_{\PP(TM)}$ and its
 bicharacteristic foliation $\Bc$ is just the null-geodesic foliation $\Lc$. This shows
 that the space $L$ of null-geodesics carries a canonical contact structure 
 $\Theta = \Theta_L$. 
 
 Explicitly, $\Theta$ can be defined as follows. Let $\gamma\in L$ be a null-geodesic,
 considered as a $1$-dimensional complex submanifold in $M$.
 At any $x\in\gamma$, the line $T_x\gamma\subset T_xM$ is isotropic, so its orthogonal
 $(T_x\gamma)^\perp$ is a hyperplane in $T_xM$ containing $T_x\gamma$.
 Now, the  contact hyperplane $\Theta_\gamma\subset  T_\gamma L$
 consists of infinitesimal displacements of $\gamma$ which, for each $x\in\gamma$,
 move $x$ inside  $(T_x\gamma)^\perp$. 
 
 \vskip .2cm
 
 We further recall that an $m$-dimensional smooth submanifold $W$ of a $2m+1$-dimensional
 contact manifold $(X,\Theta)$ is called {\em Legendrian}, if the tangent spaces of $W$ lie in $\Theta$.
 In our case $X=L$, it is clear from the above explicit definition of $\Theta$ that any subvariety
 $L_x$, see \eqref{eq:Lx}, is Legendrian.
 The main result of the ambitwistor description of conformal metrics 
 \cite{lb-3d, LB-tams, mason} can be summarized as
 follows.
 
 \begin{thm}\label{thm:ambi}
 (a) Any local Legendrian deformation of any $L_x$ inside $L$ is of the form $L_y$ for some $y$. 
 
 \vskip .2cm
 
 (b) Let $M_1, M_2$ be two civilized holomorphic conformal manifolds of dimension $n$.
 Let $x_i\in M_i$, $i=1,2$, be  two points.
 Holomorphic conformal diffeomorphisms
  $M_1\to M_2$ taking $x_1\mapsto x_2$  are in bijection with holomorphic
 contact diffeomorphisms $L(M_1) \to L(M_2)$ taking $L_{x_1}\to L_{x_2}$. \qed
 \end{thm}
 
 \begin{rem}
 In dimensions $\geq 4$ one can drop the reference to the contact structures
 thus reducing the problem to purely holomorphic geometry. That is,
 any holomorphic deformation of any $L_x$ is automatically Legendrian, and any
 holomorphic diffeomorphism $L(M_1)\to L(M_2)$ is automatically contact, see 
 \cite{LB-tams}. For special consideration of the case of dimension $3$ see
 \cite{LB-3d-old, lb-3d}. 
 Nevertheless, it seems natural to keep track of the contact structure is all dimensions,
 since it is the natural geometric structure present in the problem. 
 \end{rem}
 
 
  \section{Overcoming Liouville}
  
  \paragraph{Contact Hamiltonians and conformal Killing fields.} 
  Let $(X,\Theta)$ be a holomorphic contact manifold of dimension $2m+1$. 
  A {\em contact vector field} on $X$ is a  holomorphic 
  vector field $\xi$ preserving the distribution
  $\Theta$. That is, if we trivialize $\kap=TX/\Theta$ and view the contact form $\theta$
  as a usual 1-form, then we should have  $\Lie_\xi(\theta) = f\cdot\theta$
  for some function $f$.  It is well known \cite{ovsienko} that such a $\xi$
  is uniquely determined by the {\em contact Hamiltonian}
  \be
  H\,=\, \theta (\xi) \,\in \, H^0(X,\kap),
  \ee
 where we now view $\theta$ as a $\kap$-valued $1$-form. In this way the sheaf
 of contact vector fields is identified with the sheaf of holomorphic sections of the
 line bundle $\kap$. The Lie algebra strucure on contact vector fields
  translates to a canonical bi-differential operator (Poisson-Jacobi bracket)
   $\kap \times\kap\to\kap$. 
  
  \vskip .2cm
  
  We now specialize to $X=L(M)$ where $M$ is a civilized 
  holomorphic conformal manifold
  of dimension $n$. Theorem \ref{thm:ambi}(b) gives, as the infinitesimal version,
  the following.
  
 \begin{cor}
 We have an identification of Lie algebras
 \[
 \conf(M) \,\simeq \, H^0(L(M),\kap).  \qed
 \]
 \end{cor}
 
 \paragraph{The derived conformal algebra.} 
 The above suggests the following definition.
 
 \begin{Defi}
 Let $(M,(C_x))$ be a civilized holomorphic conformal manifold. The
 {\em derived conformal algebra} of $M$ is the dg-Lie algebra
 \[
 R\conf(M) \, := \, R\Gamma(L(M),\kap). 
 \]
 \end{Defi}
 In particular, the $1$st cohomology of this  dg-Lie algebra is $H^1(L(M),\kap)$
 which is the space of infinitesimal  deformations of $L(M)$ as a contact manifold,
 i.e., by Theorem \ref{thm:ambi},  of infinitesimal deformations of $M$ as a conformal manifold. 
 We now see that the infinite-dimensionality of the $2$-dimensional conformal group
 does not ``disappear'' in dimensions $\geq 3$, but is transformed into the 
 infinite-dimensionality of the moduli space of local conformal metrics. 
 Indeed,  symmetries and deformations are, from the point of view of derived geometry
 \cite{HS-def},
 always governed by the same algebraic structure:  an
 appropriate differential graded Lie algebra. 
 
 
 \section{The derived conformal algebra of the flat space}
 
 \paragraph{Statement of the result.}
 We now analyze the cohomology of the derived conformal algebra of the $n$-dimensional
 flat space, $n\geq 3$. We will be interested in the algebraic skeleton of the problem,
 i.e., in dealing with polynomials rather than power series. Therefore we will
 work with the algebraic variety $\AAA^n$ instead of the complex manifold $\CC^n$,
 and understand $L(\AAA^n)$ as an algebraic variety as well. So we form the
 dg-Lie algebra
 \[
 R\conf(\AAA^n)\,=\, R\Gamma(L(\AAA^n), \kap), 
 \]
considering $\kap$ as the sheaf of regular sections on the Zariski topology of $L(\AAA^n)$.  
We will identify the cohomology if this dg-Lie algebra, i.e., $H^\bullet (L(\AAA^n),\kap)$
as a module over the orthogonal group $SO(n,\CC)$. 

\vskip .2cm

More precisely, we denote by $M=\CC^n$ the standard $n$-dimensional complex vector space 
and think of $\AAA^n$ as  
``$M$ considered an an 
algebraic variety'', i.e., as the spectrum of the algebra $S^\bullet(M^*)$. 
Fixing a nondegenerate quadratic form $q\in S^2(M^*)$, we get a flat metric  on $\AAA^n$
and the variety $L(\AAA^n)$.  

\vskip .2cm

Recall the basics of representation theory of $GL(n,\CC)$,
see \cite{fulton-harris}. 
Given a sequence of integers $a=(a_1\geq \cdots \geq a_n)$
(a dominant weight for $GL(n)$), we have the {\em Schur functor}
$\Sigma^a$ from the category of $n$-dimensional $\CC$-vector spaces and their
isomorphisms to the category of finite-dimensional $\CC$-vector spaces,
with $\Sigma^a(V)$ being ``the'' space of irreducible representation of $GL(V)$
with highest weight $a$.
If all $a_i\geq 0$, we think of $a$ as a Young diagram with rows of lengths $a_1,\cdots, a_n$. 
 If $a=(a_1,\cdots, a_p, 0,\cdots, 0)$,
we write $\Sigma^{a_1,\cdots, a_p}$ for $\Sigma^a$, dropping the zeroes at the end.
We also write $1^p = ( \overbrace{1,\cdots, 1}^\text{p}, 0,\cdots, 0)$, $p\leq n$. 
Note the particular cases and properties:
\[
\begin{gathered}
\Sigma^d(V) = \Sigma^{d,0,\cdots, 0}(V) = S^p(V), \quad 
\Sigma^{1^p} (V) \,=\, \Lambda^p(V), 
\\
\Sigma^{a_1,\cdots, a_n}(V)^* \,\, \simeq\,\,  \Sigma^{a_1,\cdots, a_n}(V^*)
\,\, \simeq \, \, \Sigma^{-a_n, \cdots, -a_1}(V). 
\end{gathered}
\] 
Given two weights $a=(a_1,\cdots, a_n)$ and $b=(b_1,\cdots, b_n)$,
the decomposition of the tensor product
\[
\Sigma^a(V) \otimes \Sigma^b(V) \,\simeq \,\bigoplus_c \Sigma^c(V)^{\oplus N_{ab}^c}
\]
into irreducibles is given by the Littlewood-Richardson rule, with $N_{ab}^c$ known
as the Littlewood-Richardson coefficients. There are two important cases when 
$N_{ab}^c=1$.

\begin{exas}\label{ex:young-mult}
(a) (Horizontal Young multiplication) $c=a+b$, i.e., $c_i=a_i+b_i$. If all $a_i, b_i\geq 0$,
then $c$ is the Young diagram obtained by ``adding'' $a$ and $b$ in the horizontal
direction (row by row). The resulting projection
$\Sigma^a(V)\otimes \Sigma^b(V)\to\Sigma^{a+b}(V)$ is induced, via the Borel-Weil
theorem, by tensor multiplication of line bundles on the flag variety.

\vskip .2cm

(b) (Vertical Young multiplication). Dually, suppose that $a,b$ are nonnegative and
$c$ is the Young diagram obtained by adding $a$ and $b$ in the vertical direction,
column by column. Then $N_{ab}^c=1$ as well. The resulting projection
$y: \Sigma^a(V)\otimes \Sigma^b(V)\to\Sigma^c(V)$ can be called the
{\em vertical Young multiplication}. For instance, if $a= 1^r$,
$b= 1^s$, then
$c= 1^{r+s}$ and we get the exterior multiplication. 
We will be particularly interested in the projection
\be\label{eq:proj-y}
y_d:    S^d(V)  \otimes S^2(V)  \lra \Sigma^{d,2}(V), \quad d\geq 2. 
\ee
\end{exas}

We now specialize to $V=M^*$ where $(M,q)$ is as above and write
$SO(n) = SO(n,\CC)$ for the group of automorphisms of $(M,q)$ with determinant $1$.
Note that $M\simeq M^*$ as an $SO(n)$-module. 
The projection 
\eqref{eq:proj-y} gives an $SO(n)$-equivariant map
\[
y_{d,q} =  y(-\otimes q): S^d(M^*) \lra \Sigma^{d,2}(M^*).
\]

\begin{thm}\label{thm:conf-flat}
The dg-Lie algebra $R\conf(\AAA^n)$ has the following cohomology spaces:
\begin{itemize}
\item $H^0 \,=\, \Lambda^2(M^*)\oplus M^*\oplus M^* \oplus \CC \,=\, \Lambda^2(M^*\oplus
\CC^2) \,=\, \so(n+2)$ (the usual conformal algebra). 

\item $H^1\,=\bigoplus_{d\geq 2} \Coker (y_{d,q})$, with each $y_{d,q}$, $d\geq 2$,
being injective. 

\item $H^i=0$ for $i\geq 2$. 
\end{itemize}
\end{thm}


\paragraph{Moduli space interpretation.} We now explain why the space $H^1$ in
Theorem \ref{thm:conf-flat} can be seen as the   space of local deformations
of the conformal class of the flat metric. For this we think of the components
of  a Riemannian metric $g_{ij}(z)$
as a formal Taylor series on $\CC^n$ near $0$ and view the symmetric
algebras below as the spaces of polynomials dense in the spaces of power
series.

\vskip .2cm

The symmetric algebra $S^\bullet(M^*)$ is (after completion) the space of
formal germs of functions on $M=\CC^n$ near $0$. So the corresponding
space of germs of the metric itself is 
 the tensor product $S^2(M^*) \otimes S^\bullet(V^*)$. 
 As this is a linear space, we view it as the space of infinitesimal deformations
 of the flat metric. 
 The Pieri formula
 \cite{fulton-harris} gives
 \be\label{eq:S2xSd}
 S^2\otimes S^d \,\simeq \, S^{d+2} \oplus \Sigma^{d+1,1} \oplus \Sigma^{d,2}. 
 \ee
Let us  now quotient by changes of coordinates (understood infinitesimally, as vector fields).
The space of vector fields (understood in the same sense as above) is
$M\otimes S^\bullet(M^*)$. We identify $M$ with $M^*$ as a $SO(n)$-module. 
Again, the Pieri formula gives
\be\label{eq:MxSd}
M^*\otimes S^d(M^*) \,\simeq \, S^{d+1}(M^*) \oplus \, \Sigma^{d,1}(M^*).
\ee
 So the ``moduli space" of metrics modulo coordinate changes
 has, as the
 tangent space at the trivial metric, the 
 result of subtracting the contributions from \eqref{eq:MxSd} for all $d$
 from the contributions from \eqref{eq:S2xSd} for all $d$, which gives 
   $\bigoplus_{d\geq 2} \Sigma^{d,2}(M^*)$. For instance, the lowest summand here,
   $\Sigma^{2,2}(M^*)$, is precisely the space of all possible values of
   the Riemann curvature tensor at the origin. 
   
   Further, let us look at the effect of passing to conformal classes, i.e., quotienting
   by multiplication by scalar functions, on the tangent space to the moduli space.
   The space of functions is $S^\bullet(V^*)$.  So taking the cokernel of the map
   \[
   y_q: S^{\geq 2}(V) \lra \bigoplus_{d\geq 2} \Sigma^{d,2}(V)
   \]
 has the effect of passing to the tangent space of the moduli space of conformal
 classes.

 \begin{rem}\label{rem:2d}
 Finally, it is instructive to compare the situation with the $2$-dimensional case
 when we have an infinite-dimensional conformal algebra in homological degree $0$.
 The difference is that for $\dim(M)=2$ the map
 \[
 y_{d,q}: S^d(M^*) \lra \Sigma^{d,2}(M^*) =  S^{d-2}(M^*) \otimes
  \Lambda^2(M^*)^{\otimes 2}
 \]
 is surjective, not injective. The kernel of $y_{d,q}$
 has dimension $2$, it is the space of traceless symmetric tensors in $2$ variables.
 So in each degree we have  two basis vectors contributing to 
 the kernel.  This matches the identification
 \[
 \conf(\AAA^2) \,=\, \CC[z] \del_z \oplus \CC[\ol z] \del_{\ol z}. 
 \]
\end{rem}


\section {Proof of Theorem  \ref{thm:conf-flat}.} 
\paragraph{ Identifying the bundle $\kap$. } 
We first identify the line bundle $\kap$, the target of the contact form,
starting from the compact flat case.
That is, let $V=\CC^{n+2}$ with a non-degenerate scalar product 
$\langle -,- \rangle$ and let $Q^n\subset \PP(V)=\PP^{n+1}$
be the quadric of null-directions. Then $L(Q^n)$ is the isotropic Grassmannian 
$G^\is(2,V)\subset G(2,V)$. We denote by $S$ the tautological rank $2$ bundle
on both $G(2,V)$ and $G^\is(2,V)$ and put $\Oc(1)=\Lambda^2(S^*)$.

\begin{lem}\label{lem:kap=O(1)}
The line bundle $\kap_{G^\is(2,V)}$ is identified with $\Oc(1)$. 
\end{lem}

\noindent{\sl Proof:} Let $E\subset V$ be a $2$-dimensional isotropic subspace
and $[E]\in G^\is(2,V)$ be the corresponding point. Then it is standard that
\[
T_{[E]} G(2,V)\,\simeq \, \Hom(E, V/E).
\]
Inside this, $T_{[E]}G^\is(2,V)$ consists of linear maps $f: E\to V/E$
such that 
\be\label{eq:f(e),e}
\langle f(e),e\rangle =0 \quad \text{for any } e\in E.   
\ee
(Since $E$ is isotropic, $\langle f(e),e\rangle$ is well defined.) This is
a codimension $3$ subspace in $\Hom(E,V/E)$. Further, the contact hyperplane
$\Theta_E\subset T_{[E]} G^\is(2,V)$ is   $\Hom(E, E^\perp/E)$ (a codimension $4$
subspace in  $\Hom(E,V/E)$),
see the general discussion in \S \ref{sec:ambi}B. 
 Given $f$ satisfying \eqref {eq:f(e),e}, we have
 \[
 \langle f(e_1), e_2\rangle \,=\, -\langle f(e_2), e_1\rangle, \quad 
 \text{for any } e_1, e_2\in E.  
 \]
 Therefore the expression $ \langle f(e_1), e_2\rangle$ is a linear map
 $\Lambda^2(E)\to\CC$. Vanishing of this map means that $f: E\to E^\perp/E$,
 i.e., $f\in\Theta_E$.
 This gives an identification of vector spaces
 \[
 T_{[E]} G^\is(2,V)/\Theta_E \,\simeq \, \Lambda^2(E^*),
 \]
 and so an identification of line bundles $\kap\simeq \Lambda^2(S^*) = \Oc(1)$. \qed
 
 \vskip .2cm
 
 We now pass from $L(Q^n)$ to the Zariski open part $L(\AAA^n)$ which is, by 
 \eqref {eq:L(Cn)}, the total space of an algebraic vector bundle whose
 projection we denote by $\pi$:
 \[
 L(\AAA^n) \,=\, \Tot \bigl( T_{\PP^{n-1}}(-1)|_{Q^{n-2}}\bigr) 
\buildrel \pi\over\lra Q^{n-2}. 
 \]
 Let us write for short 
 \[
 Q:=Q^{n-2},\quad 
 G\, := \, (T_{\PP^{n-1}}(-1))^* = \Omega^1_{\PP^{n-1}}(1).
  \]
 In a more algebro-geometric language
 the identification of $L(\AAA^n)$ with the total space reads:
 \[
 L(\AAA^n) \,=\,\Spec \,\bigoplus_{d\geq 0} S^d(G|_Q). 
 \]
Lemma \ref{lem:kap=O(1)} implies that
\[
\kap_{L(\AAA^n)} \,\simeq \, \pi^* \Oc_Q(1),
\] 
 and therefore
 \[
 H^i(L(\AAA^n,\kap) \,=\,\bigoplus_d H^i(Q, S^d(G) (1)|_Q). 
 \]
 
 \paragraph{Cohomology on $\PP(M)$ using Borel-Weil-Bott.}
 We invoke the short exact sequence of sheaves on $\PP^{n-1}=\PP(M)$
 \be\label{eq:onPand Q}
 0\to S^d(G)(-1) \buildrel \cdot q\over\lra S^d(G)(1) \lra S^d(G)(1)|_Q \to 0
 \ee
  and analyze first the cohomology of $S^d(G)(\pm 1)$ on $\PP(M)$. 
  
  \begin{lem}\label{lem:onP(M)}
  On $\PP(M)$, 
  \begin{itemize}
  \item[$(-1)$] The sheaf $S^d(G)(-1)$ has $H^1 =  S^{d-1}(M^*)$ (understood as $0$
  for $d=0$) and no other cohomology. 
  
  \item[$(+1)$] 
  \begin{itemize}
  \item The sheaf $S^0(G)(1)$ has $H^0=M^*$ and no other cohomology.
  \item The sheaf $S^1(G)(1)$ has $H^0=\Lambda^2(M^*)$ and no other cohomology.
  \item The sheaf $S^2(G)(1)$ has no cohomology.
  \item The sheaf $S^d(G)(1)$, $d\geq 3$, has $H^1=\Sigma^{d-1,2}(M^*)$ and no other cohomology. 
  \end{itemize}
  \end{itemize}
  \end{lem}
  
  \noindent{\sl Proof:} We use the Borel-Weil-Bott theorem for flag varieties,
  see \cite{manin}, Ch.1, \S 2 for a treatment convenient for us.

 Let $F=F(M)$ be the space of complete flags
 \[
 M_1\subset M_2 \subset \cdots \subset M_n=M=\CC^n,\quad \dim(M_i)=i, 
 \]
 with the natural projection $p: F\to \PP(M) = \{M_1\subset M\}$. 
 We denote by $M_i$ the tautological bundle on $F$ of rank $i$. 
  To a weight $a=(a_1, \cdots, a_n)\in\ZZ^n$ (not necessarily dominant)
  we associate the line bundle
  \[
  \Oc_F(a)\,=\, (M/M_{n-1})^{\otimes a_1} \otimes (M_{n-1}/M_{n-2})^{\otimes a_2}
  \otimes\cdots\otimes (M_2/M_1)^{\otimes a_{n-1}} \otimes M_1^{a_n}
  \]
 on $F$.  As mentioned in  Example \ref{ex:civilized}(a), $G^*=T_{\PP(M)}(-1)$
 is the universal quotient bundle whose fiber at $M_1\subset M$ is $M/M_1$.
 This implies that
 \[
 S^d(G) \,=\, p_* \, \Oc_F (0,\cdots, 0, -d, 0).
 \]
 Indeed, taking the space of sections of  the line bundle $(M_2/M_1)^{\otimes (-d)} = \Oc_{\PP(M/M_1)}(d)$
 on the projective space $\PP(M/M_1)$ or, equivalently, of the pullback of this
 line bundle to the full flag variety of $M/M_1$, gives $S^d(M/M_1)^*$. 
 This implies that for any $b\in\ZZ$
 \[
 S^d(G)(b) \,=\, p_* \Oc_F(0, \cdots, 0, -d,b),
 \]
 and so 
 \be\label{eq:SdGb}
 H^\bullet (\PP(M), S^d(G)(b)) \,=\, H^\bullet (F, \Oc_F(0,\cdots, 0, -d, b)). 
 \ee
 We now recall the procedure of finding $H^\bullet(F, \Oc_F(a))$ for $a\in\ZZ^n$
 given by Bott's theorem. That is, if $\lambda = (\lambda_1\geq \cdots \geq \lambda_n)$
 is a dominant weight, $w\in S_n$ is a permutation of length $\ell(w)$ and
 $\rho=(n, n-1, \cdots, 1)$, then
 \[
 H^i\bigl(F, \Oc_F(w(\lambda+\rho)-\rho)\bigr) \,=\,
 \begin{cases}
 \Sigma^\lambda(M),& \text{ if } i=\ell(w),
 \\
 0,&\text{ otherwise.}
 \end{cases}
 \]
 Thus to find $H^\bullet(F,\Oc_F(a))$ we need to represent $a=w(\lambda+\rho)-\rho$
 with $\lambda$ dominant.  If such a representation is impossible, i.e., if $a+\rho$
 has repetitions, then $\Oc_F(a)$ has no cohomology. 
 
 \vskip .2cm
 
 After these preparations, let us establish part $(+1)$ of Lemma \ref {lem:onP(M)}.
 From \eqref{eq:SdGb} we see that we need to find $H^\bullet(F, \Oc_F(a))$,
 where
 \[
 a=(0,\cdots, 0, -d,1), \quad \text{so} \quad a+\rho = (n, n-1, \cdots, 3, 2-d,2). 
 \]
 For $d=0$ we have a repetition so no cohomology. For $d\geq 1$,
 a single elementary transposition (length $1$) takes $a+\rho$ to
 $(n, n-1,\cdots, 3, 2, 2-d)$, then subtracting $\rho$ we get $(0,\cdots, 0, 1-d)$.
 So in this case the only non-trivial cohomology is
 \[
 H^1(F, \Oc_F(a)) \,=\,\Sigma^{0, \cdots, 0, 1-d}(M) \,=\, S^{d-1}(M^*)
 \]
 as claimed. 
 
 \vskip .2cm
 
 Let us now establish part $(-1)$ of Lemma \ref {lem:onP(M)}. We have
 \[
  a=(0,\cdots, 0, -d,-1), \quad \text{so} \quad a+\rho = (n, n-1, \cdots, 3, 2-d,0). 
 \]
 Now,
 \begin{itemize}
 \item If $d=0$, then $a$ is dominant so we have only
 \[
 H^0(F, \Oc_F(a))\,=\, \Sigma^{0,\cdots, 0, -1}(M)\,=\, M^*.
 \]
 \item If $d=1$, then $a$ is still dominant, so we have only 
  \[
 H^0(F, \Oc_F(a))\,=\, \Sigma^{0,\cdots, 1, -1}(M)\,=\, \Lambda^2(M^*).
 \]
 \item If $d=2$, we get $a+\rho=(\cdots, 3,0,0)$, a repetition so no cohomology.
 
 \item If $d\geq 3$, then $a+\rho=(\cdots, 3, 2-d,0)$ which is  ordered, by an
 elementary transposition, to $(\cdots, 3,0,2-d)$. Subtracting $\rho$, we get
 $(0,\cdots, 0, -2, 1-d)$, so   the only cohomology is
 \[
 H^1(F, \Oc_F(a)) \,=\, \Sigma^{0,\cdots, 0, -2, 1-d}(M) \,=\, \Sigma^{d-1,2}(M^*). 
 \]
 \end{itemize}
 Lemma \ref {lem:onP(M)} is proved. \qed
 
 \paragraph{Cohomology on $Q$.} 
 We now finish the proof of Theorem \ref{thm:conf-flat}. Let us display the cohomology
 (known from Lemma \ref {lem:onP(M)})
 of the first two sheaves $S^d(G)(\pm 1)$ in  
 \eqref{eq:onPand Q} in a table (Fig. \ref{fig:1}), under these sheaves.  Under the third sheaf, $S^d(G)(1)|_Q$,
 let us write
 the conclusion about its cohomology preceded by the sign ``$\Rightarrow$''.
 We note that in the last row,  the map
  $H^1(\PP(M), S^d(G)(-1))\to H^1(\PP(M), S^d(G)(1))$
 induced by multiplication with $q$, is proportional to $y_{d-1,q}$.
This follows by  invariance, by letting $q\in S^2(M^*)$ vary and
using the fact (Example \ref{ex:young-mult}(b)) that
\[
\dim \, \Hom_{GL(M)} \bigl( S^{d-1}(M^*) \otimes S^2(M^*), \, \Sigma^{d-1,2}(M^*)\bigr) \,=\, 1. 
\]
The fact that the coefficient of proportionality is non-zero, is implied by the next lemma.

\begin{lem}\label{lem:injective}
The map
\[
\wt q: H^1(\PP(M), S^d(G)(-1)) \lra H^1(\PP(M), S^d(G)(1))
\]
induced by multiplication with $q$, is injective. 
\end{lem} 
This lemma, together with the table in Fig. \ref{fig:1},  establish Theorem  \ref{thm:conf-flat}.
 
 \begin{figure}[H]
  \[
 \xymatrix{
 0\ar[r]& S^d(G)(-1) \ar[r]^{\cdot q}& S^d(G)(1) \ar[r]& S^d(G)(1)|_Q \ar[r]& 0
 \\
 d=0&  H^\bullet=0  & H^0=M^* & \Rightarrow \, H^0=M^* &
 \\
 d=1& H^1=\CC & H^0=\Lambda^2(M^*) & \Rightarrow \, H^0=\CC\oplus \Lambda^2(M^*) &
 \\
 d=2 & H^1=M^* & H^\bullet=0 & \Rightarrow \, H^0=M^* &
 \\
 d\geq 3 & H^1=S^{d-1}(M^*) \ar[r]^{c\cdot y_{d-1,q}} & H^1=\Sigma^{d-1,2}(M^*) & 
 \Rightarrow H^1 = \Coker(y_{d-1,q}).
 }
 \]
 \caption{Calculating cohomology on $Q\subset \PP(M)$.}\label{fig:1}
 \end{figure}
 
 \noindent{\sl Proof of Lemma \ref{lem:injective}:}  Let $\varpi: F\to G(2,M)$ be
 the projection. We denote by $M_2$ the tautological rank $2$ bundle on $G(2,M)$. 
 If $E\subset M$ is a $2$-dimensional subspace and $[E]\in G(2,M)$ is the corresponding
 point, then $\varpi^{-1}(E) = F(M/E)\times \PP(E)$. By applying the Borel-Weil-Bott
 theorem to the fibers of $\varpi$, we find that $\wt q$ is  identified with the morphism
 \[
 S^{d-1}(M^*) = H^0\bigl(G(2,M), S^{d-1}(M_2)^*\bigr) \lra H^0\bigl(G(2,M),\Sigma^{d-1,2}(M_2^*)\bigr)
 = \Sigma^{d-1,2}(M^*)
 \]
 induced by the morphism of vector bundles on $G(2,M)$
 \[
 y_{d-1, q|_{M-2}}:  S^{d-1}(M_2)^*) \lra ,\Sigma^{d-1,2}(M_2^*)
 \]
 which, on each fiber, i.e., on each $E\subset M$ as above, is
 the morphism
 \[
 y_{d-1, q|_E}: S^{d-1}(E^*)\lra \Sigma^{d-1,2}(E^*)
 \]
 corresponding to the $2$-dimensional space $E$ and the quadratic form $q|_E$.
 This morphism has been discussed in Remark \ref{rem:2d}, and its kernel
 is the subspace  in $S^{d-1}(E^*)$ formed by polynomials harmonic (traceless)
 with respect to $q|_E$. So we are reduced to the following fact.
 
 \begin{lem}
 Let $M$ be a complex vector space of dimension $\geq 3$ and $q\in S^2(M^*)$ be
 a non-degenerate quadratic form. If $f\in S^{d-1}(M^*)$, $d\geq 3$,  is such that for any
 $2$-dimensional subspace $E\subset M$,  the restriction
 $f|_E$ is harmonic with respect to  $q|_E$, then $f=0$. 
 \end{lem}

 \noindent{\sl Proof:} We can assume $q$ to come from a positive definite quadratic
 form on a real form $M_\RR=\RR^n$  of $M$. Then it is enough to prove the lemma 
 under the assumptions that $f$ is a real homogeneous polynomial of 
 degree $d-1$ on $\RR^n$
 and the restriction of $f$ to any real subspace $E$ is harmonic with respect to $q|_E$. 
 If $E$ is a $2$-dimensional real space with a positive definite quadratic form,
 then we can use Euclidean geometry in $E$. In particular, a harmonic
 polynomial homogeneous of degree $m$ is, in polar coordinates $(R,\phi)$
 a linear combination of $R^m \cos(m\phi)$ and $R^m \sin(m\phi)$, and therefore
  it is invariant under Euclidean rotations by $2\pi/m$ in $E$.
 So our assumptions on $f: \RR^n \to\RR$ imply that the restriction of $f$
 to any $2$-plane $E\subset \RR^n $ is invariant under rotations by $2\pi/(d-1)$
 in this plane. If $d\geq 4$, this implies that $f(x)$ depends only on the radius
 $\|x\| = q(x)^{1/2}$, so by homogeneity $f(x)=\on{const}\cdot \|x\| ^{d-1}$,
 which contradicts the above trigonometric shape of $f|_E$, so $f=0$. 
 
 In the remaining case $d=3$, our $f$ is a quadratic form. The condition that $f|_E$
 is harmonic with respect  to $q|_E$ means, in the classical terminology, that
 $f|_E$ and $q|_E$ are ``anharmonic'', i.e., that the sum of the two
 eigenvalues of $f_E$ with respect to $q|_E$ is $0$. Since $q$ is assumed
 to be positive definite, this implies that  each $f|_E$ must
 always  have signature $(+,-)$, if non-degenerate and must be zero, if degenerate.
  This is impossible unless $f=0$. 
  \qed


M.K.: Kavli IPMU, 5-1-5 Kashiwanoha, Kashiwa, Chiba, 277-8583 Japan, \\
{\tt mikhail.kapranov@protonmail.com}

\ed